\documentclass{article}

\usepackage{ifxetex,ifluatex}
\newif\ifxetexorluatex
\ifxetex
  \xetexorluatextrue
\else
  \ifluatex
    \xetexorluatextrue
  \else
    \xetexorluatexfalse
  \fi
\fi

\ifxetexorluatex
  \usepackage{fontspec}
\else
  \usepackage[T1]{fontenc}
  \usepackage[utf8]{inputenc}
  \usepackage{lmodern}
\fi

\usepackage{amsfonts,amsmath,graphicx}
\newtheorem{thm}{Theorem}

\title{A parametrization of $8\times8$ magic squares of squares
 through octonionic multiplication}
\author{Ísabel Pirsic}
\date{}

\begin{document}
\maketitle

In 1771 Euler posed what he called a `problema curiosum' \cite{eulQuat}: to find
a $4\times 4$ magic square of squares,
i.e., an integer matrix $M$ such that $MM^\top=M^\top M=c*I,\ c\in\mathbb Z$,
where also the sums of the squares of the diagonal entries should add up to $c$.

He gave explicit solutions, e.g.,
\begin{center}\begin{tabular}{|c|c|c|c|}
    \hline
    $68$  & $-29$ & $41$  & $-37$ \\ \hline
    $-17$ & $31$  & $79$  & $32$  \\ \hline
    $59$  & $28$  & $-23$ & $61$  \\ \hline
    $-11$ & $-77$ & $8$   & $49$  \\
    \hline
  \end{tabular}\end{center}
but also the following interesting parametrisation:
{{\begin{center}\begin{tabular}{|c|c|c|c|}
        \hline
        +a p+b q+      & +a r$-$b s      & $-$a s$-$b r      & +a q$-$b p      \\
        +c r+d s      & $-$c p+d q &      +c q+d p &      +c s$-$d r \\ \hline
        $-$a q+b p+      & +a s+b r      & +a r$-$b s      & +a p+b q      \\
         +c s$-$d r    &  +c q+d p    &  +c p$-$d q    &  $-$c r$-$d s    \\ \hline
        +a r+b s$-$      & $-$a p+b q      & +a q+b p      & +a s$-$b r      \\
         $-$c p$-$d q    &  $-$c r+d s    &  +c s+d r    &  $-$c q+d p    \\ \hline
        $-$a s+b r$-$     & $-$a q$-$b p      & $-$a p+b q      & +a r+b s      \\
         $-$c q+d p    &  +c s+d r    &  +c r$-$d s    &  +c p+d q    \\
        \hline
      \end{tabular}\end{center} }}
where, to fulfill the diagonal conditions, two equations have to be
satisfied, namely
\begin{equation*} \label{conds}
  p r+q s\quad=\quad0,\qquad \frac{a}c =
  \frac{-d(pq+rs)-b(ps+qr)}{b(pq + rs) + d(ps + qr)}.
\end{equation*}

\begin{figure}
{\centering
\includegraphics[width=0.75\textwidth]{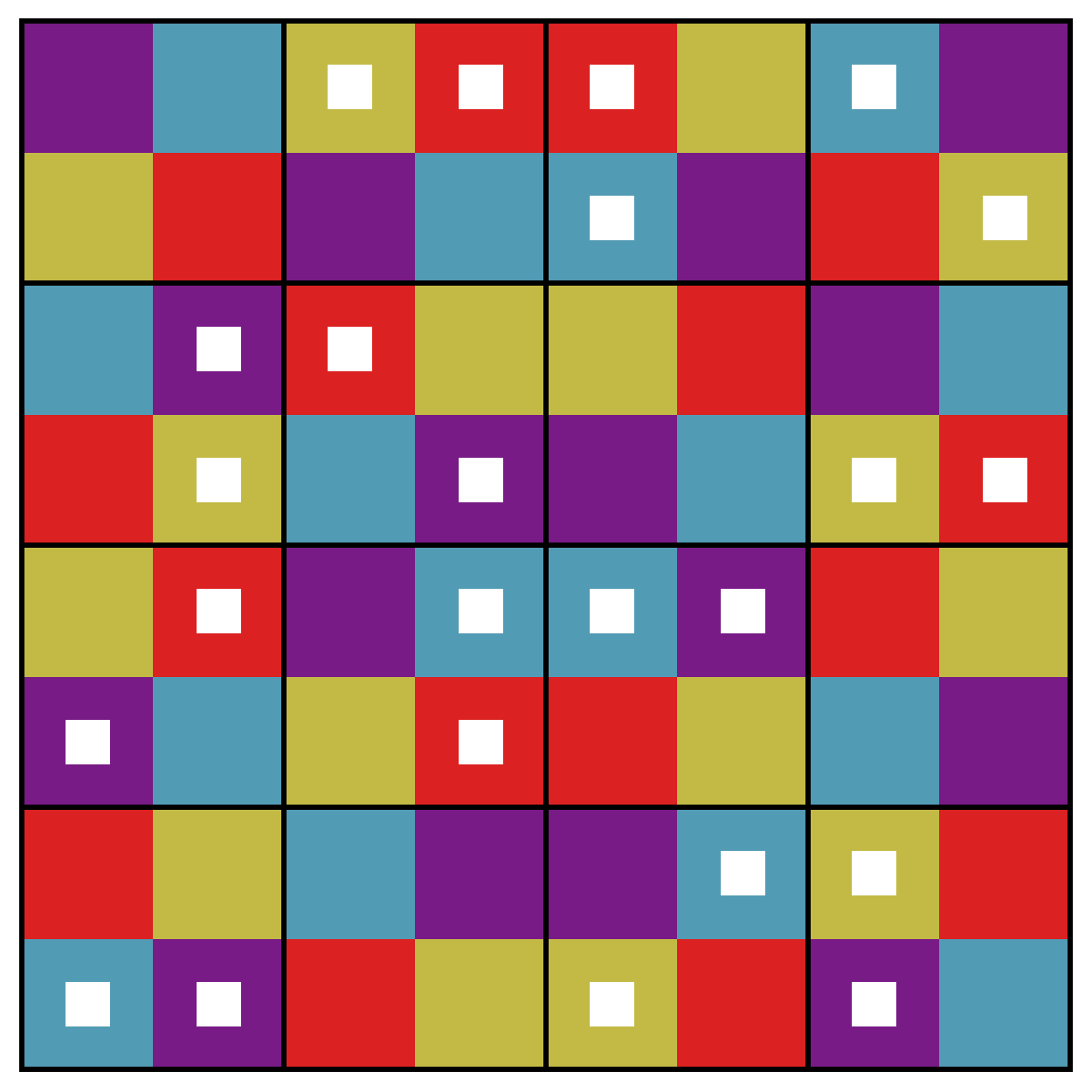}
\caption{\label{quatMat}The quaternionic magic square of squares}}
\end{figure}

It is not known exactly, how Euler arrived at the solution, though we can
take some guesses. For instance, the
visualization presented in Figure \ref{quatMat}
hints at a strong connection to Latin squares, famously another of
Euler's research objects \cite{eulsq}.
-- The figure is to be read like this:
For each $2\times2$ subsquare, the colors signify the first factor of a term
like $a p$, the place inside the subsquare the second and the white squares
negated terms.

Only in the early 1900s, Hurwitz gave a deeper explanation for Euler's matrix.
In the last chapter of his `Vorlesungen über die Zahlentheorie der Quaternionen'
\cite{hurwQuat}, he effectively proves that $SO(4,\mathbb R)$ can be
modelled by combined left and right
multiplication of quaternions, i.e., for any $\gamma\in SO(4)$ there exist
quaternions $q_L,q_R$ of norm 1, such that
\[
  \gamma\big((a,b,c,d) \big) = (p,q,r,s) \quad\iff\quad
  q_L*(a+b i+c j+d k)*q_R = p+qi+rj+sk ,
\]
where $i,j,k$ are the quaternionic basis elements. (We will
identify 4-vectors and quaternions.)
Therefore, if we initially have a orthonormal basis
$\{\boldsymbol e_1,
\boldsymbol e_2, \boldsymbol e_3, \boldsymbol e_4\}= \{1,i,j,k\}$, transformation
by multiplication to
$$A*\boldsymbol e_i*P= (a+bi+cj+dk)*\boldsymbol e_i*(p+qi+rj+sk),$$ with
not necessarily unit quaternions $A,P$ will give an at least orthogonal basis
again. And this is exactly what happens in Euler's matrix, save some sign
changes and column permutation;  the rows are given by the
components of $A*\boldsymbol e_i*P$. The column permutation is
necessary to make the diagonal conditions easier to attain, sign changes are
probably just for aesthetic reasons.

Euler did not know yet about quaternions, but he had found the four-squares-identity
\[
\begin{array}{c}
{\displaystyle (a_{1}^{2}+a_{2}^{2}+a_{3}^{2}+a_{4}^{2})(b_{1}^{2}+b_{2}^{2}+b_{3}^{2}+b_{4}^{2})=}\\
{\displaystyle (a_{1}b_{1}-a_{2}b_{2}-a_{3}b_{3}-a_{4}b_{4})^{2}+}
{\displaystyle (a_{1}b_{2}+a_{2}b_{1}+a_{3}b_{4}-a_{4}b_{3})^{2}+}\\
{\displaystyle (a_{1}b_{3}-a_{2}b_{4}+a_{3}b_{1}+a_{4}b_{2})^{2}+}
{\displaystyle (a_{1}b_{4}+a_{2}b_{3}-a_{3}b_{2}+a_{4}b_{1})^{2}}
\end{array}
\]
which secretly is also founded on quaternionic multiplication. It may indeed
be seen as the multiplicativity of the quaternionic norm, i.e., $N(A)N(B)=N(AB)$.
This is relevant in the magic square context, since the norm of the rows is thus
constant and equal to \[ N(A\boldsymbol e_i P) = N(A)N(\boldsymbol e_i)N(P)=
N(A)N(P) .\]

Now arose the question, if a similar thing can be done for the octonions,
the infamously non-associative Cayley-Dickson extension of the quaternions.
There, it is known  (see, e.g., \cite{conwOct}) that $SO(8)$ is generated by at most
seven multiplications, either all left or all right, but conversely
any single multiplication with a
unit octonion is also associated to a  transformation in $SO(8)$.

\begin{figure}
{\centering
\includegraphics[width=0.9\textwidth]{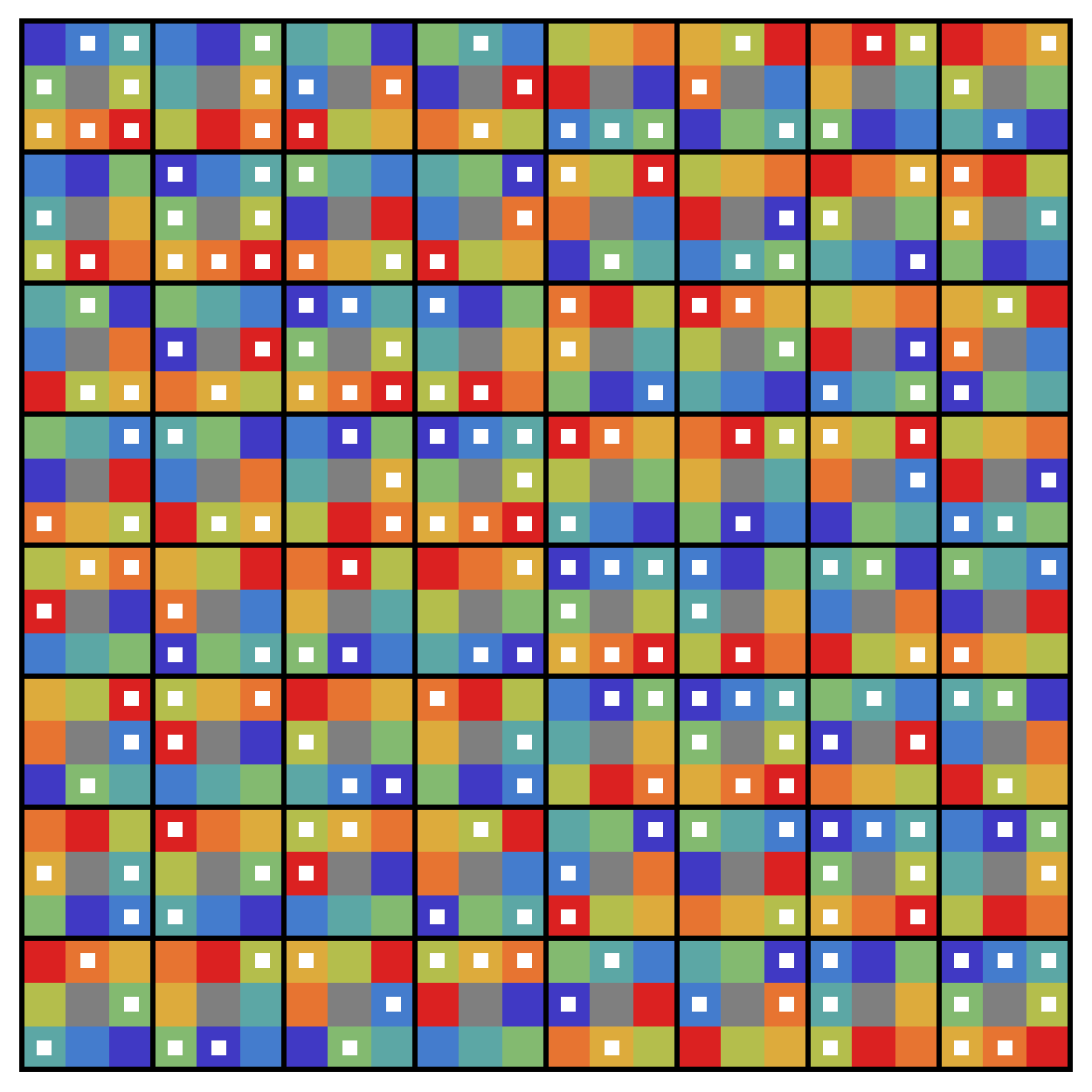}
\caption{\label{octMat}The octonionic semimagic square of squares}}
\end{figure}

So, in order to obtain a magic square of size 8 with reasonably small norm,
the same approach as shown by Hurwitz may be made, i.e., to build a matrix
from the components of $A*(\boldsymbol e_i*P),$ where $\boldsymbol e_i$ runs
through the octonionic basis elements. (Multiplying with just
one octonion would not give a `magic' matrix, i.e., where all entries are distinct;
so two is the minimum. Similarly, performing two multiplications on the same side
does not produce enough distinct elements.) -- To
give this matrix in detail in text form
is not particularly enlightning, we will instead present its visualization analogous
to the $4\times 4$ case (the grey fields are just a way to complete
the square) in Figure \ref{octMat}.

Some properties are apparent in this visualization, e.g., the various
Sudoku-like combinatorial conditions that are fulfilled in the rows and columns
of subsquares, and also the quadrants,
a central symmetry regarding subsquares, as well as, on a finer
scale, the white squares representing the negative terms. Intriguing is also the fact
that offside the diagonal all choices of 3 elements from 8 are present.

This matrix is then indeed found to be row-orthogonal, as can also be verified using
a computer algebra system; leading to the only `theorem' in this note:
\begin{thm}
 Let
 \begin{align*}
  A=&\ (a\boldsymbol e_1+b\boldsymbol e_2+c\boldsymbol e_3+d\boldsymbol e_4+
  e\boldsymbol e_5+f\boldsymbol e_6+g\boldsymbol e_7+h \boldsymbol e_8),\\
  P=&\ (p\boldsymbol e_1+q\boldsymbol e_2+r\boldsymbol e_3+s\boldsymbol e_4+
  t\boldsymbol e_5+u\boldsymbol e_6+v\boldsymbol e_7+w \boldsymbol e_8),\\
  A, P\in&\
   \mathbb O[a,b,\dots,h,\,p,q,\dots,w]
  \end{align*}
  be symbolic octonions, (i.e., $a,b,\dots,h,\,p,q,\dots,w$ are variables, and
 $\boldsymbol e_1=1,\,\boldsymbol e_i$ are orthogonal basis
 elements of $\mathbb O$) and let
 $M\in \mathbb R[a,\dots,h,p,\dots,w]^{8\times8}$ be given by the $i$-th row vector
 defined as the components of $A(\boldsymbol e_i P)$.

 Then the matrix $M$ constitutes a symbolic semi-magic square of squares
 with entries in $\mathbb Z[a,\dots,h,p,\dots,w]$,
 with orthogonal rows and constant $N(A)N(P)$.

 Consequently, supposing the matrix is evaluated at some $A_0,P_0\in \mathbb Q^8$,
 such that all entries are distinct integers, this integer matrix
 constitutes an (integral)
 semi-magic square of squares with orthogonal rows and constant $N(A_0)N(P_0)$.
\end{thm}
\textsc{Proof\ }
 The assertions on $M$ follow from the fact that multiplication with
 unit octonions models $SO(8)$ (see \cite{conwOct}), hence right multiplication with
 $P$ and subsequent left multiplication with $A$ preserves the initial orthogonality,
 while multiplying all lengths equally. Also, all choices of orthogonal
 basis elements act as a permutation on each other, up to sign change, so that
 the coefficients of the matrix entries are all $\pm1$.
 That all entries are different
 can be verified by calculation, e.g., using Sage \cite{sage}. Figure \ref{octMat}
 also visualizes this in that no two $3\times3$-subsquares have the same pattern
 combination of colors and white squares. The rest of the
 assertions then follows.\hfill$\Box$

-- Having attained a parametrisation the
next step is to find numbers such that the entries are indeed different and the
constant is low. By plugging in concrete values for the
$A,P$, e.g., this matrix with constant 9476 could be found:
\begin{center} \begin{tabular}{|r|r|r|r|r|r|r|r|}\hline
     43& 16& -19& 8& -22& 47& 38& -53 \\ \hline
     -30& 11& 30& 5& 25& -32& 75& -16 \\ \hline
     9& -4& -7& -52& -46& -57& -6& -35 \\ \hline
     -8& -67& 48& 21& -5& 10& -17& -42 \\ \hline
     54& -11& 14& 49& -31& -36& 17& 36 \\ \hline
     44& 41& 60& -21& 33& -2& -25& -10 \\ \hline
     7& -26& 29& -54& -24& 37& 32& 45 \\ \hline
     -41& 46& 35& 22& -60& 15& -12& 1\\ \hline
\end{tabular} \end{center}
which only lacks the diagonal properties required for a fully magic matrix;
this type is called semi-magic.

The diagonal properties in Euler's matrix arose from choosing such a permutation
of rows that in the diagonal all terms were present. Then calculating the
sum of squares of the diagonal elements and subtracting the magic constant leads
very naturally to the side conditions. -- This does not work with the $8\times8$
matrix, at least not in an obvious straightforward manner, so that
for now we cannot ascertain whether
or not a fully magic matrix exists. A small scale randomized
computer search using parameters within the range $-32$ to $32$ did
not produce any results.

It is not difficult to find at least parameters such that not only all entries
but also all their squares (or absolute values) are distinct, one of low constant
43617 is given by
$$A=(-2, -3, 7, -1, 2, -11, 2, -5), \quad P=(-7, 4, -4, -9, 1, 2, -5, 3).$$
When presenting the matrix in terms of the parameters, it is dependent on the
specific implementation of the algebra. We chose the Cayley-Dickson extension
of the quaternion algebra as given by the computer algebra system Sage \cite{sage}.
(The according scripts are available from the author by request.)

The parameters for the 9476 matrix exhibit an interesting detail, in that one of
them has half integers for coordinates:
 \[  A= (8, -2, -4, 8, -4, -1, -5, -4),\quad
   P=\frac12(5,7,-1 ,-3 ,-7 ,1 ,7 ,1) .\]
It is already the case in the quaternions that the appropriate integers to
perform number theory contain half integers, this is even more true for the
octonions (see a concise account in \cite{conwOct}).
The half integers may be placed only in specific 4-subsets of
the coordinates. It seemed appropriate to choose at least one factor of this
form to achieve a low constant.

The next step by the Cayley-Dickson construction would be the 16-dimensional
sedenions; however, the Hurwitz theorem states that there are no further
normed real division algebras beyond dimension 8.
In particular the multiplicativity
of the norm does not hold any longer; consequently the analog matrix does not
have orthogonal rows. Similar constructions might also be carried out in
any other composition algebra, though a different definition
of orthogonality would apply, i.e., instead of $MM^\top$ one would
consider something like $MQM^\top$ for a matrix $Q$.


\end{document}